\documentstyle[11pt]{article}
\begin{document}
\setlength{\textheight}{574pt}
\setlength{\textwidth}{432pt}
\setlength{\oddsidemargin}{18pt}
\setlength{\topmargin}{14pt}
\setlength{\evensidemargin}{18pt}
\newtheorem{theorem}{Theorem}[section]
\newtheorem{lemma}{Lemma}[section]
\newtheorem{corollary}{Corollary}[section]
\newtheorem{conjecture}{Conjecture}
\newtheorem{remark}{Remark}[section]
\newtheorem{definition}{Definition}[section]
\newtheorem{problem}{Problem}
\newtheorem{example}{Example}
\newtheorem{proposition}{Proposition}[section]
\title{{\bf INVARIANCE OF PLURIGENERA OF VARIETIES WITH NONNEGATIVE KODAIRA DIMENSIONS}}
\date{November 28, 2000}
\author{Hajime Tsuji}
\maketitle
\begin{abstract}
We prove the invariance of plurigenera under smooth projective 
deformations of projective varieties with 
nonnegative Kodaira dimensions. 
MSC32J25
\end{abstract}
\tableofcontents
\section{Introduction}
Let $X$ be a smooth projective variety and let $K_{X}$ 
be the canonical bundle of $X$.
The canonical ring 
\[
R(X,K_{X}):= \oplus_{m\geq 0}H^{0}(X,{\cal O}_{X}(mK_{X}))
\]
is a basic birational invariant of $X$.  
For every positive integer $m$, 
the $m$-th plurigenus $P_{m}(X)$ is defined by 
\[
P_{m}(X) : = \dim H^{0}(X,{\cal O}_{X}(mK_{X})). 
\]
It is believed that $P_{m}(X)$ is invariant under 
smooth projective deformation. 

Recently Y.-T. Siu (\cite{si}) proved that $P_{m}(X)$ 
is invariant under smooth projective deformation, 
if the all the fibers are of general type. 
This result has been slightly generalized by \cite{ka,nak}.

The main idea of Siu is the comparison of the AZD on 
a special fiber and one on the total space of the deformation 
by using an induction on multiplier ideal sheaves. 

In this paper, we consider 
the case of a smooth projective deformation 
such that all the fibers are smooth projective 
varieties with nonnegative Kodaira dimensions. 
\begin{theorem}
Let $\pi : X \longrightarrow \Delta$ be a 
smooth projective family of smooth projective varieties 
with nonnegative Kodaira dimensions 
over the open unit disk in the complex plane with 
center $0$. 
Then for every positive integer $m$,  
the $m$-th plurigenus $P_{m}(X_{t}) (X_{t} : = \pi^{-1}(t))$
 is independent of $t\in \Delta$. 
\end{theorem}
We also prove the following theorem. 
\begin{theorem}
Let $\pi : X \longrightarrow \Delta$ be a 
smooth projective family of smooth projective varieties 
with pseudoeffective canonical bundles over the open 
unit disk $\Delta$ in the complex plane with center $0$. 
Then the numerical Kodaira dimension $\kappa_{num}(X_{t})$ 
is independent of $t\in \Delta$. 
\end{theorem}
Let us explain the idea of the proof of Theorem 1.1. 
Let $\pi : X\longrightarrow \Delta$ be the family 
as in Theorem 1.1.
The heart of the proof is to construct a singular hermitian metric 
$h$ on $K_{X}$ such that $\Theta_{h}$ is semipositive and 
$h\mid_{X_{t}}$ is an AZD of $K_{X_{t}}$ for every $t\in \Delta$. 
Then by the $L^{2}$-extension theorem (\cite[p.200, Theorem]{o-t}),
we may easily deduce that $P_{m}(X_{t})$ is independent of 
$t \in \Delta$.  

Let $A$ be an ample line bundle on $X$. 
Let $m$ be an arbitrary positive integer. 
By the induction as in \cite{si}, it is easy to see that 
for a canonical AZD $h_{m}$ of $A + mK_{X}$, 
$h_{m}\mid_{X_{t}} (t\in \Delta )$  is an AZD of 
${\cal O}_{X_{t}}(A+mK_{X})$ for every $t\in \Delta$.
But this is much weaker than what we want. 

The key point of the proof is to embed   $X$ into a 
family of strongly pseudoconvex manifolds  $B$ as a divisor with the negative normal bundle
$-A$. 
Then by the Bergman construction of an AZD of $K_{B}+ X$ on $B$,
we can construct the  desired AZD of $K_{X}$.  

This technique can be used to prove the abundance of canonical divisor of 
a projective manifold with nonnegative Kodaira dimension. 
This will be treated in the forthcoming paper. 
 
\section{Preliminaries}
\subsection{Multiplier ideal sheaves}
In this subsection $L$ will denote a holomorphic line bundle on a complex manifold $M$. 
\begin{definition}
A  singular hermitian metric $h$ on $L$ is given by
\[
h = e^{-\varphi}\cdot h_{0},
\]
where $h_{0}$ is a $C^{\infty}$-hermitian metric on $L$ and 
$\varphi\in L^{1}_{loc}(M)$ is an arbitrary uppersemicontinuous function on $M$.
We call $\varphi$ a  weight function of $h$.
\end{definition}
The curvature current $\Theta_{h}$ of the singular hermitian line
bundle $(L,h)$ is defined by
\[
\Theta_{h} := \Theta_{h_{0}} + \sqrt{-1}\partial\bar{\partial}\varphi ,
\]
where $\partial\bar{\partial}$ is taken in the sense of a current.
The $L^{2}$-sheaf ${\cal L}^{2}(L,h)$ of the singular hermitian
line bundle $(L,h)$ is defined by
\[
{\cal L}^{2}(L,h) := \{ \sigma\in\Gamma (U,{\cal O}_{M}(L))\mid 
\, h(\sigma ,\sigma )\in L^{1}_{loc}(U)\} ,
\]
where $U$ runs over the  open subsets of $M$.
In this case there exists an ideal sheaf ${\cal I}(h)$ such that
\[
{\cal L}^{2}(L,h) = {\cal O}_{M}(L)\otimes {\cal I}(h)
\]
holds.  We call ${\cal I}(h)$ the {\bf multiplier ideal sheaf} of $(L,h)$.
If we write $h$ as 
\[
h = e^{-\varphi}\cdot h_{0},
\]
where $h_{0}$ is a $C^{\infty}$ hermitian metric on $L$ and 
$\varphi\in L^{1}_{loc}(M)$ is the weight function, we see that
\[
{\cal I}(h) = {\cal L}^{2}({\cal O}_{M},e^{-\varphi})
\]
holds.
For $\varphi\in L^{1}_{loc}(M)$ we define the multiplier ideal sheaf of $\varphi$ by 
\[
{\cal I}(\varphi ) := {\cal L}^{2}({\cal O}_{M},e^{-\varphi}).
\]
Also we define 
\[
{\cal I}_{\infty}(h) := {\cal L}^{\infty}({\cal O}_{M},e^{-\varphi})
\]
and call it the $L^{\infty}$-{\bf multiplier ideal sheaf} of $(L,h)$.

\begin{definition}
$L$ is said to be pseudoeffective, if there exists 
a singular hermitian metric $h$ on $L$ such that 
the curvature current 
$\Theta_{h}$ is a closed positive current.

Also a singular hermitian line bundle $(L,h)$ is said to be pseudoeffective, 
if the curvature current $\Theta_{h}$ is a closed positive current.
\end{definition}
It is easy to see that a line bundle $L$ on a smooth projective 
manifold $M$ is 
pseudoeffective,
if and only if for an ample line bundle $H$ on $M$,  
$L + \epsilon H$ is {\bf Q}-effective (or big) 
for every positive rational number $\epsilon$ (cf. \cite{d}).

The following theorem is fundamental in the applications 
of multiplier ideal sheaves. 
\begin{theorem}(Nadel's vanishing theorem \cite[p.561]{n})
Let $(L,h)$ be a singular hermitian line bundle on a compact K\"{a}hler
manifold $M$ and let $\omega$ be a K\"{a}hler form on $M$.
Suppose that $\Theta_{h}$ is strictly positive, i.e., there exists
a positive constant $\varepsilon$ such that
\[
\Theta_{h} \geq \varepsilon\omega
\]
holds.
Then ${\cal I}(h)$ is a coherent sheaf of ${\cal O}_{M}$ ideal 
and for every $q\geq 1$
\[
H^{q}(M,{\cal O}_{M}(K_{M}+L)\otimes{\cal I}(h)) = 0
\]
holds.
\end{theorem}

We note that the multiplier ideal sheaf of a singular hermitian {\bf R}-line 
bundle (i.e., a real power of a line bundle) is well defined because the multiplier ideal sheaf is defined 
in terms of the weight function.
Sometimes it is useful to consider the 
following variant of multiplier ideal sheaves. 
\begin{definition}
Let $h_{L}$ be a singular hermitian metric on a line bundle $L$.
Suppose that the curvature of $h_{L}$ is a positive current on $X$.
We set 
\[
\bar{\cal I}(h_{L}) 
:= \lim_{\varepsilon\downarrow 0}{\cal I}(h_{L}^{1+\varepsilon})
\]
and call it   the closure of ${\cal I}(h_{L})$. 
\end{definition} 
As you see later, the closure of a multiplier ideal  sheaf is easier to handle 
than the original multiplier ideal sheaf in some respect. 

Next we shall consider the restriction of singular hermitian 
line bundles to subvarieties. 
\begin{definition}
Let $h$ be a singular hermitian metric on $L$ given by 
\[
h = e^{-\varphi}\cdot h_{0},
\]
where $h_{0}$ is a $C^{\infty}$-hermitian metric on $L$ and 
$\varphi\in L^{1}_{loc}(M)$ is an uppersemicontinuous function.
Here $L^{1}_{loc}(M)$ denotes the set of locally integrable functions 
(not the set of classes of almost everywhere equal locally 
integrable functions on $M$). 

For a subvariety $V$ of $M$, we say that the restriction 
$h\mid_{V}$ is well defined, if 
$\varphi$ is not identically $-\infty$ on $V$. 
\end{definition}
Let $(L,h)$,$h_{0}$,$V$, $\varphi$ be as in Definition 2.4.
Suppose that the curvature current $\Theta_{h}$ is bounded 
from below by some $C^{\infty}$-(1,1)-form. 
Then $\varphi$ is an almost plurisubharmonic 
function, i.e. locally a sum of a plurisubharmonic function and a $C^{\infty}$-function.
Let $\pi :\tilde{V} \longrightarrow V$ be an arbitrary 
resolution of $V$. 
Then  $\pi^{*}(\varphi\mid_{V})$ is locally integrable on $\tilde{V}$, 
since $\varphi$ is almost plurisubharmonic. 
Hence 
\[
\pi^{*}(\Theta_{h}\mid_{V}) := \Theta_{\pi^{*}h_{0}\mid_{V}}
+ \sqrt{-1}\partial\bar{\partial}\pi^{*}(\varphi\mid_{V})
\]
is well defined.

\begin{definition}
Let $\varphi$ be a plurisubharmonic function on a unit 
open polydisk  $\Delta^{n}$ with center $O$.
We define the Lelong number of $\varphi$ at $O$ by 
\[
\nu (\varphi ,O) := \liminf_{x\rightarrow O}\frac{\varphi (x)}{\log \mid x\mid},
\]
where $\mid x\mid  = (\sum\mid x_{i}\mid^{2})^{1/2}$.
Let $T$ be a closed positive $(1,1)$-current on 
a unit open polydisk $\Delta^{n}$.
Then by $\partial\bar{\partial}$-Poincar\'{e} lemma
there exists a plurisubharmonic function  $\phi$ 
on $\Delta^{n}$ such that
\[
T = \frac{\sqrt{-1}}{\pi}\partial\bar{\partial}\phi .
\]
We define the Lelong number $\nu (T,O)$ at $O$ by
\[
\nu (T,O) := \nu (\phi ,O).
\]
It is easy to see that $\nu (T,O)$ is independent of the choice of
$\phi$ and local coordinates around $O$.
For an analytic subset $V$ of a complex manifold $X$, we set 
\[
\nu (T,V) = \inf_{x\in V}\nu (T,x).
\]
\end{definition}
\begin{remark} More generally 
the Lelong number is defined for a closed positive
$(k,k)$-current on a complex manifold.
\end{remark}

\begin{theorem}(\cite[p.53, Main Theorem]{s})
Let $T$ be a closed positive $(k,k)$-current on a complex manifold
$M$.
Then for every $c > 0$
\[
\{ x\in M\mid \nu (T,x)\geq c\}
\]
is a subvariety of codimension $\geq k$
in $M$.
\end{theorem}

The following lemma shows a rough relationship between 
the Lelong number of $\nu(\Theta_{h},x)$ at $x\in X$ and the stalk of the multiplier
ideal sheaf ${\cal I}(h)_{x}$ at $x$. 

\begin{lemma}(\cite[p.284, Lemma 7]{b}\cite{b2},\cite[p.85, Lemma 5.3]{s})
Let $\varphi$ be a plurisubharmonic function on 
the open unit polydisk $\Delta^{n}$ with center $O$.
Suppose that $e^{-\varphi}$ is not locally integrable 
around $O$.
Then we have that
\[
\nu (\varphi ,O)\geq 2
\]
holds.
And if
\[
\nu (\varphi ,O) > 2n
\]
holds, then $e^{-\varphi}$ is not locally integrable around $O$.
\end{lemma}
Let $(L,h)$ be a pseudoeffective singular hermitian line bundle 
on a complex manifold $M$. 
The {\bf closure} $\bar{\cal I}(h)$ of the multiplier ideal 
sheaf ${\cal I}(h)$ can be analysed in terms of Lelong numbers
in the following way.  
We note that $\bar{\cal I}(h)$ is coherent ideal sheaf on $M$ by Theorem 2.1.

In the case of $\dim M = 1$, we can compute $\bar{\cal I}(h)$ in terms of the Lelong number $\nu (\Theta_{h},x) (x\in M)$.  
In fact in this case $\bar{\cal I}(h)$ is locally free and 
\[
\bar{\cal I}(h) = {\cal O}_{M}(-\sum_{x\in M}[\nu (\Theta_{h},x)]x)
\]
holds by Lemma 2.1, because $2  = 2\dim M$. 

In the case of $\dim M \geq 2$,
let $f : N \longrightarrow M$ be a modification such that 
$f^{*}\bar{\cal I}(h)$ is locally free. 
If we take $f$ properly, we may assume that there exists a divisor $F = \sum_{i} F_{i}$ with normal crossings on $Y$ such that 
\[
K_{N} = f^{*}K_{M} + \sum_{i} a_{i}F_{i}
\]
and 
\[
\bar{\cal I}(h) = f_{*}{\cal O}_{N}(-\sum_{i} b_{i}F_{i})
\]
hold on $Y$ for some nonnegative integers $\{ a_{i}\}$ and $\{ b_{i}\}$. 
Let $y \in F_{i}- \sum_{j\neq i}F_{j}$ and let 
$(U,z_{1},\ldots ,z_{n})$ be a local corrdinate neighbourhood  of $y$ which is biholomorphic to the open unit disk $\Delta^{n}$ with center $O$ in $\mbox{\bf C}^{n}(n = \dim M)$ and 
\[
U \cap F_{i} = \{ p\in U\mid z_{1}(p) = 0\}
\]
holds.
For $q\in \Delta^{n-1}$, we  set $\Delta (q):= \{ p\in U\mid 
(z_{2}(p),\ldots ,z_{n}(p)) = q\}$.
Then considering the family of the restriction $\{ \Theta_{h}\mid_{\Delta (q)}\}$ for very general $q\in \Delta^{n-1}$, by Lemma 2.1, we see that   
\[
b_{i} = \max\{ [\nu (f^{*}\Theta_{h},F_{i}) - a_{i}],0\}
\]
holds for every $i$. 
In this way $\bar{\cal I}(h)$ is determined by the {\bf Lelong numbers} of 
the curvature current on some modification. 
This is not the case, unless we take the closure as in the following example. 
\begin{example}
Let $h_{P}$ be a singular hermitian metric on the trivial line bundle on 
the open unit polydisk $\Delta$ with center $O$ in {\bf C} 
defined by
\[
h_{P} = \frac{\parallel\,\cdot\,\parallel^{2}}{\mid z\mid^{2}(\log \mid z\mid  )^{2}}.
\]
Then $\nu (\Theta_{h_{P}},0) = 1$ holds. 
But ${\cal I}(h_{P}) = {\cal O}_{\Delta}$ holds. 
On the other hand $\bar{\cal I}(h_{P}) = {\cal M}_{0}$ holds, 
where ${\cal M}_{0}$ is the ideal sheaf of $0\in \Delta$.  
\end{example}

\subsection{Restriction of multiplier ideal sheaves to divisors}

Let $(L,h)$ be a pseudoeffective singular hermitian line bundle on 
a smooth projective variety $X$. 
Let $D$ be a smooth divisor on $X$. 
We set  
\[
v_{m}(D) = \mbox{mult}_{D}\mbox{Spec}({\cal O}_{X}/{\cal I}(h^{m}))
\]
and 
\[
\tilde{\cal I}_{D}(h^{m}) = {\cal O}_{D}(v_{m}(D)D)\otimes {\cal I}(h^{m}).
\] 
Then $\tilde{\cal I}_{D}(h^{m})$ is an ideal sheaf on $D$ (it is torsion free, 
since $D$ is smooth). 

Let $x\in D$ be an arbitrary point of $D$ and let 
$(U,z_{1},\ldots ,z_{n}) (n:= \dim X)$ be a local coordinate neighbourhood of 
$x$ which is 
biholomorphic to the unit open polydisk $\Delta^{n}$ with center $O$ in 
$\mbox{\bf C}^{n}$ and 
\[
U \cap D = \{ p\in U\mid z_{1}(p) = 0\}
\] 
holds. 
For $q\in \Delta^{n-1}$, we  set $\Delta (q):= \{ p\in U\mid 
(z_{2}(p),\ldots ,z_{n}(p)) = q\}$.
Then considering the family of the restriction $\{ \Theta_{h}\mid_{\Delta (q)}\}$ for very general $q\in \Delta^{n-1}$, by Lemma 2.1, we see that  
\[
m\cdot\nu (\Theta_{h},D) -1 \leq v_{m}(D) 
\leq m\cdot\nu (\Theta_{h},D)
\]
holds. 

We define the ideal sheaves $\sqrt[m]{\tilde{\cal I}_{D}(h^{m})}$ 
on $D$ by 
\[
\sqrt[m]{\tilde{\cal I}_{D}(h^{m})}_{x}
:= \cup {\cal I}(\frac{1}{m}(\sigma))_{x} (x\in D),
\]
where $\sigma$ runs all the germs of $\tilde{\cal I}_{D}(h^{m})_{x}$.
And we set 
\[
{\cal I}_{D}(h):= \cap_{m\geq 1}\sqrt[m]{\tilde{\cal I}_{D}(h^{m})}
\]
and call it {\bf the multipler ideal of $h$ on $D$}.
Also we set 
\[
\bar{\cal I}_{D}(h):= \lim_{\varepsilon\downarrow 0}
{\cal I}_{D}(h^{1+\varepsilon}).
\]
The following theorem is crucial in our proof of Theorem 1.1.
\begin{theorem}(\cite[Theorem 2.8]{tu4}) 
Let $(L,h)$ be a singular hermitian line bundle on a smooth projective 
variety $X$. 
Suppose that $\Theta_{h}$ is bounded from below by some 
negative multiple of a $C^{\infty}$-K\"{a}hler form on $X$. 
Let $D$ be a smooth divisor on $X$. 
If $h\mid_{D}$ is well defined, then 
\[
\bar{\cal I}_{D}(h) = \bar{\cal I}(h\mid_{D})
\]
holds. 
\end{theorem}

\subsection{Analytic Zariski decompositions}
In this subsection we shall introduce the notion of analytic Zariski decompositions. 
By using analytic Zariski decompositions, we can handle  big line bundles
like  nef and big line bundles.
\begin{definition}
Let $M$ be a compact complex manifold and let $L$ be a holomorphic line bundle
on $M$.  A singular hermitian metric $h$ on $L$ is said to be 
an analytic Zariski decomposition, if the followings hold.
\begin{enumerate}
\item $\Theta_{h}$ is a closed positive current,
\item for every $m\geq 0$, the natural inclusion
\[
H^{0}(M,{\cal O}_{M}(mL)\otimes{\cal I}(h^{m}))\rightarrow
H^{0}(M,{\cal O}_{M}(mL))
\]
is an isomorphim.
\end{enumerate}
\end{definition}
\begin{remark} If an AZD exists on a line bundle $L$ on a smooth projective
variety $M$, $L$ is pseudoeffective by the condition 1 above.
\end{remark}

\begin{theorem}(\cite{tu,tu2})
 Let $L$ be a big line  bundle on a smooth projective variety
$M$.  Then $L$ has an AZD. 
\end{theorem}
As for the existence for general pseudoeffective line bundles, 
now we have the following theorem.
\begin{theorem}(\cite{d-p-s})
Let $X$ be a smooth projective variety and let $L$ be a pseudoeffective 
line bundle on $X$.  Then $L$ has an AZD.
\end{theorem}
{\bf Proof of Theorem 2.5}. Let  $h_{0}$ be a fixed $C^{\infty}$-hermitian metric on $L$.
Let $E$ be the set of singular hermitian metric on $L$ defined by
\[
E = \{ h ; h : \mbox{lowersemicontinuous singular hermitian metric on $L$}, 
\]
\[
\hspace{70mm}\Theta_{h}\,
\mbox{is positive}, \frac{h}{h_{0}}\geq 1 \}.
\]
Since $L$ is pseudoeffective, $E$ is nonempty.
We set 
\[
h_{L} = h_{0}\cdot\inf_{h\in E}\frac{h}{h_{0}},
\]
where the infimum is taken pointwise. 
The supremum of a family of plurisubharmonic functions 
uniformly bounded from above is known to be again plurisubharmonic, 
if we modify the supremum on a set of measure $0$(i.e., if we take the uppersemicontinuous envelope) by the following theorem of P. Lelong.

\begin{theorem}(\cite[p.26, Theorem 5]{l})
Let $\{\varphi_{t}\}_{t\in T}$ be a family of plurisubharmonic functions  
on a domain $\Omega$ 
which is uniformly bounded from above on every compact subset of $\Omega$.
Then $\psi = \sup_{t\in T}\varphi_{t}$ has a minimum 
uppersemicontinuous majorant $\psi^{*}$ which is plurisubharmonic.
\end{theorem}
\begin{remark} In the above theorem the equality 
$\psi = \psi^{*}$ holds outside of a set of measure $0$(cf.\cite[p.29]{l}). 
\end{remark}

By Theorem 2.6 we see that $h_{L}$ is also a 
singular hermitian metric on $L$ with $\Theta_{h}\geq 0$.
Suppose that there exists a nontrivial section 
$\sigma\in \Gamma (X,{\cal O}_{X}(mL))$ for some $m$ (otherwise the 
second condition in Definition 3.1 is empty).
We note that  
\[
\frac{1}{\mid\sigma\mid^{\frac{2}{m}}} 
\]
gives the weihgt of a singular hermitian metric on $L$ with curvature 
$2\pi m^{-1}(\sigma )$, where $(\sigma )$ is the current of integration
along the zero set of $\sigma$. 
By the construction we see that there exists a positive constant 
$c$ such that  
\[
\frac{h_{0}}{\mid\sigma\mid^{\frac{2}{m}}} \geq c\cdot h_{L}
\]
holds. 
Hence
\[
\sigma \in H^{0}(X,{\cal O}_{X}(mL)\otimes{\cal I}_{\infty}(h_{L}^{m}))
\]
holds.  
Hence in praticular
\[
\sigma \in H^{0}(X,{\cal O}_{X}(mL)\otimes{\cal I}(h_{L}^{m}))
\]
holds.  
 This means that $h_{L}$ is an AZD of $L$. 
\vspace{10mm} {\bf Q.E.D.} 
\begin{remark}
By the above proof we have that for the AZD $h_{L}$ constructed 
as above
\[
H^{0}(X,{\cal O}_{X}(mL)\otimes{\cal I}_{\infty}(h_{L}^{m}))
\simeq 
H^{0}(X,{\cal O}_{X}(mL))
\]
holds for every $m$. 
\end{remark}
The following proposition implies that the multiplier ideal sheaves 
of $h_{L}^{m}(m\geq 1)$ constructed in the proof of
 Theorem 2.5 are independent of 
the choice of the $C^{\infty}$-hermitian metric $h_{0}$.
The proof is trivial.  Hence we omit it. 
\begin{proposition}
$h_{0},h_{0}^{\prime}$ be two $C^{\infty}$-hermitian metrics 
on a pseudoeffective line bundle $L$ on a smooth projective 
variety $X$. 
Let $h_{L},h^{\prime}_{L}$ be the AZD's constructed as in the 
proof of Theorem 2.5 associated with $h_{0},h_{0}^{\prime}$ 
respectively. 
Then 
\[
(\min_{x\in X}\frac{h_{0}}{h_{0}^{\prime}}(x))\cdot h_{L}^{\prime}
     \leq    h_{L} \leq
 (\max_{x\in X}\frac{h_{0}}{h_{0}^{\prime}}(x))\cdot h_{L}^{\prime}
\]
hold.
In particular 
\[
{\cal I}(h_{L}^{m}) = {\cal I}((h_{L}^{\prime})^{m})
\]
holds for every $m\geq 1$.
\end{proposition}
We call the AZD constructed as in the proof of Theorem 2.5  {\bf a canonical 
AZD} of $L$. 
Proposition 2.1 implies that the multiplier ideal sheaves associated with 
the multiples of the canonical AZD are independent of the choice of 
the canonical AZD. 

\subsection{Numerical Kodaira dimension}
The numerical Kodaira dimension was first introduced by Y. Kawamata in \cite{ka}, for nef line bundles on a projective algebraic varieties.
Here we give another definition of the numerical Kodaira dimension 
for any line bundles on projective algebraic varieties. 
This definition coincides with the former definition when the 
line bundle is nef. 
\begin{definition}
Let $L$ be a line bundle on a smooth projective variety
$X$.  Let $A$ be an  ample line bundle on $X$. 
We define the numerical Kodaira dimension $\kappa_{num}(L)$ of $L$ by
\[
\kappa_{num}(X,L) :=
\sup_{\ell\geq 1}(\limsup_{m\rightarrow\infty}\frac{\log \dim H^{0}(X,{\cal O}_{X}(\ell A+mL))}{\log m}).
\]
And we define the numerical Kodaira dimension $\kappa_{num}(X)$ of $X$ by 
\[
\kappa_{num}(X) := \kappa (X,K_{X}).
\]
\end{definition}

\begin{remark}
It is clear that $\kappa_{num}(X,L)\geq 0$, if and only if $L$ is pseudoeffective. 
If $L$ is nef, then 
\[
\kappa_{num}(X,L) = \sup\{ k\mid c_{1}^{k}(L)\,\,\mbox{ is not numerically trivial}  \}
\]
holds.  The righthandside is the definition of the numerical Kodaira dimension 
defined in \cite{ka}. 
\end{remark}
It is clear that 
\[
\kappa_{num}(X,L) \geq \kappa (X,L)
\]
holds.

\section{Proof of Theorem 1.2} 
Let $\pi : X\longrightarrow \Delta$ be the 
family as in Theorem 1.2. 
Let $A$ be an ample line bundle on $X$ 
Let $h_{m,0}$ be a canonical AZD of 
$mK_{X_{0}}+A$ and let $h_{m}$ be a canonical AZD of 
$mK_{X} +A$. 
Then we have the following lemma.
\begin{lemma}
For every positive integer $\ell$
\[
\bar{\cal I}(h_{m,0}^{\ell}) = 
\bar{\cal I}(h_{m}^{\ell}\mid_{X_{0}})
\]
holds. 
\end{lemma}
{\bf Proof of Lemma 3.1.} 
We  prove this lemma by induction on $m$. 
If $m =0$, then the assertion is clear. 
Suppose that the assertion is settled for $m-1 (m\geq 1)$. 
Let $h_{m-1}$ be a canonical AZD of $(m-1)K_{X}+ A$.
Let $\xi_{m-1}$ be the function such that 
\[
h_{m-1} = e^{-\xi_{m-1}}\cdot H_{m-1}
\]
holds, where $H_{m-1}$ is a $C^{\infty}$-hermitian metric 
on $(m-1)K_{X}+ A$. 
Let us fix a $C^{\infty}$-hemitian metric $H_{m}$ on 
$mL +A$ 
Let  $\varphi_{0}$,$\varphi$  be functions on 
$X_{0}$ defined by
\[
h_{m,0} = e^{-\varphi_{0}}\cdot H_{m}
\]
and 
\[
h_{m}\mid _{X_{0}} = e^{-\varphi}\cdot H_{m}. 
\]
Let $L_{m}$ denotes the line bundle $mK_{X}+A$. 
Let $E$ be an effective {\bf Q}-divisor on $X$ such that 
$L_{m} - E$ is ample. 
Let $\psi$ be a function on $X$ such that 
\[
e^{-\psi}\cdot H_{m}
\]
is a $C^{\infty}$ hermitian metric on $L_{m} - E$. 

There exists a positive integer $a$ such that  
\[
M = a(L_{m}- E)
\]
is Cartier and for any pseudoeffective singular hemitian line bundle 
$(L,h_{L})$ on $X$, ${\cal O}_{X}(M + L)\otimes{\cal I}(h_{L})$ 
is globally generated.
The existence of such $a$ follows 
from \cite[p.664, Proposition 1]{si}.

From now on we consider all the functions on $X_{0}$. 
The following lemma is clear. 
\begin{lemma}
We may assume that 
\[
\psi \leq \varphi \leq \varphi_{0} \leq  \xi_{m-1} 
\]
hold, if we adjust them by adding  real constants. 
\end{lemma}
The following lemma is similar to \cite[p.670, Proposition 5]{si}.
\begin{lemma}
Let $0 < \varepsilon << 1$ be a sufficiently small positive number such that 
$e^{-\varepsilon\psi}$ is locally integrable every where on $X$ 
and $e^{-\varepsilon\psi}\mid_{X_{0}}$ is locally integrable
everywhere on $X_{0}$.
Then  
\[
{\cal I}((\ell -\varepsilon)\varphi_{0} + (a+\varepsilon )\psi )
\subset 
{\cal I}((\ell -1 + a-\varepsilon )\varphi +\varepsilon\psi + \xi_{m-1})
\]
holds for every positive integer $\ell$. 
\end{lemma}
{\bf Proof}. 
The proof here is essentially same as in \cite[p.670, Proposition 5]{si}. 
We prove the lemma by induction on $\ell$. 

If $\ell = 1$, 
\begin{eqnarray*}
(1-\varepsilon )\varphi_{0} + (a+\varepsilon )\psi 
 & \leq (1-\varepsilon )\varphi_{0} + a\varphi + \varepsilon\psi \\
 &\leq  (a-\varepsilon )\varphi + \varepsilon\psi + \xi_{m-1}
\end{eqnarray*}
hold.

Suppose that the induction step $\ell$ has been settled. 
Let $\sigma$ be an element of 
\[
H^{0}(X_{0},{\cal O}_{X_{0}}((\ell +a+1)L_{m})\otimes
{\cal I}((\ell+1 -\varepsilon)\varphi_{0} + (a+\varepsilon )\psi ).
\]
Then by the induction assumption, we have that
\[
H^{0}(X_{0},{\cal O}_{X_{0}}((\ell +a+1)L_{m})\otimes
{\cal I}((\ell+1 -\varepsilon)\varphi_{0} + (a+\varepsilon )\psi ))
\]
\[
\subset  
H^{0}(X_{0},{\cal O}_{X_{0}}((\ell +a+1)L_{m})\otimes
{\cal I}((\ell-\varepsilon)\varphi_{0} + (a+\varepsilon )\psi )) 
\]
\[
\subset  
H^{0}(X_{0},{\cal O}_{X_{0}}((\ell +a+1)L_{m})\otimes
{\cal I}((\ell -1 + a-\varepsilon )\varphi +\varepsilon\psi + \xi_{m-1})) \
\]
\[
= 
H^{0}(X_{0},{\cal O}_{X_{0}}(L_{m}+ (\ell +a)L_{m})\otimes
{\cal I}((\ell -1 + a-\varepsilon )\varphi +\varepsilon\psi + \xi_{m-1})) 
\]
hold. 
We note that the singular hermitian metric 
\[
h_{m-1} = e^{-\xi_{m-1}}\cdot H_{m-1}
\]
on $(m-1)K_{X}+A$ has semipositive curvature in the sense of a current. 
Then by the definitions  of $L_{m}$ and $\xi_{m-1}$,  Nadel's vanishing theorem 
(Theorem 2.1) implies  that 
\[
H^{1}(X,{\cal O}_{X}(L_{m}-X_{0}+ (\ell +a)L_{m})\otimes
{\cal I}((\ell -1 + a-\varepsilon )\varphi +\varepsilon\psi + \xi_{m-1})) 
= 0
\]
holds. 
Hence $\sigma$ is the restriction of an element of 
\[
H^{0}(X,{\cal O}_{X}(L_{m}+ (\ell +a)L_{m})\otimes
{\cal I}((\ell -1 + a-\varepsilon )\varphi +\varepsilon\psi + \xi_{m-1})).
\]
This implies that 
\[
\sigma \in  H^{0}(X_{0},{\cal O}_{X_{0}}((\ell +1  +a)L_{m})\otimes
{\cal I}((\ell+ a-\varepsilon )\varphi +\varepsilon\psi + \xi_{m-1}))
\]
holds by the definition of $\varphi$, i.e., by the definition of 
an AZD.  
Since 
\[
H^{0}(X_{0},{\cal O}_{X_{0}}((\ell +a+1)L_{m})\otimes
{\cal I}((\ell+1 -\varepsilon)\varphi_{0} + (a+\varepsilon )\psi )).
\]
is globally generated by the construction of $A$, we conclude that 
\[
{\cal I}((\ell+1 -\varepsilon)\varphi_{0} + (a+\varepsilon )\psi )
\subset 
{\cal I}((\ell + a-\varepsilon )\varphi +\varepsilon\psi + \xi_{m-1})
\]
holds. 
This completes the proof of Lemma 3.3.  
 {\bf Q.E.D.} \vspace{5mm} \\

Let 
\[
f : Y_{0}\rightarrow X_{0}
\]
be any modification.
By  Lemma 3.3 and Lemma 2.1,  we see that 
\[
\nu (f^{*}((\ell -\varepsilon )\varphi_{0}+ (a+\varepsilon )\psi ) ) + n
\geq \nu (f^{*}((\ell -1+a-\varepsilon ) \varphi
+\varepsilon\psi + \xi_{m-1})) 
\]
holds, where $n = \dim X_{0}$. 
Dividing the both sides by $\ell$  and letting $\ell$ tend to infinity, 
we have that 
\[
\nu (f^{*}\varphi_{0}) \geq \nu (f^{*}\varphi )
\]
holds. 
On the other hand, we note that  
\[
\nu (f^{*}\varphi_{0}) \leq \nu (f^{*}\varphi )
\]
holds by Lemma 3.2. 
Hence we see that for every $k\geq 0$, 
\[
\bar{\cal I}(k\varphi_{0}) = \bar{\cal I}(k\varphi )
\]
holds on $X_{0}$, 
since the closure of a multiplier ideal sheaf is determined 
by the Lelong numbers of pullback of the curvature currents 
on some modification of the space as in Section 2.1. 
By the definitions of $\varphi$ and $\varphi_{0}$,
this completes the proof of Lemma 3.1. 
{\bf Q.E.D.} \vspace{5mm} \\

Now we shall prove Theorem 1.2. 
By Lemma 3.1 and the $L^{2}$-extension theorem (\cite[p.200, Theorem]{o-t}),
we see that 
\[
H^{0}(X,{\cal O}_{X}(mK_{X}+2A)\otimes \bar{\cal I}(h_{m}))
\rightarrow 
H^{0}(X_{0},{\cal O}_{X_{0}}(mK_{X}+2A)\otimes \bar{\cal I}(h_{m,0}))
\]
is surjective for every $m \geq 0$. 
Since $A$ can be replaced by its any positive multiple, 
we see that 
\[
\kappa_{num}(X_{t}) \leq \kappa_{num}(X_{t})
\]
for every very general $t\in\Delta$. 
On the other hand by the uppersemicontinuity of 
the cohomology groups, we see that the opposite inequality 
holds.
Hence we see that 
$\kappa_{num}(X_{t})$ is independent of $t\in \Delta$.
This completes the proof of Theorem 1.2. 
\section{Proof of Theorem 1.1}
Let $\pi : X\longrightarrow \Delta$ be the 
family as in Theorem 1.1. 
Let $h$ be a canonical AZD of $K_{X}$ and let 
$h_{0}$ be  a canonical AZD of $K_{X_{0}}$. 

\subsection{Logarithmic reduction} 

Let $A$ be an ample line bundle on $X$ and let 
$h_{A}$ be a $C^{\infty}$-hermitian metric on $X$ 
with strictly positive curvature.
We may and do assume that 
for every pseudoeffective singular hermitian line bundle 
$(F,h_{F})$, 
${\cal O}_{X}(F+A)\otimes {\cal I}(h_{F})$ is 
globally generated (cf. \cite[Proposition 1.1]{si}). 
We set $L := {\cal O}_{X}(-A)$ and let 
$p : L \longrightarrow X$ be the bundle projection. 
Let $B$ be the unit disk bundle associated with 
the hermitian line bundle $(L,h_{A}^{-1})$. 
Then $B$ is a family of strongly pseudoconvex manifolds over $\Delta$. 
We denote the bundle projection $B\longrightarrow X$ 
again by $p$ and we shall identify $X$ with the 
zero section of $L$.  
Let $h_{B}$ be a $C^{\infty}$-hermitian metric on 
$K_{B}+X$.
Let $dV_{L}$ a $C^{\infty}$-volume form on $L$ and 
let $dV_{B}$ be the restriction of $dV_{L}$ to $B$.
We define (the diagonal part of) the $m$-th Bergman kernel $K_{m}(B)$
of $(K_{B}+ X,h_{B})$ by 
\[
K_{m}(B) = \sum_{i=1}^{\infty}\mid\phi_{i}^{(m)}\mid^{2},
\]
where $\{\phi_{i}^{(m)}\}_{i=1}^{\infty}$ is a complete orthonormal basis
of the Hilbert space ${\cal H}_{m}$ of $L^{2}$-holomorphic sections 
of $m(K_{B}+L)$ with respect to $h_{B}^{m}$ and $dV_{B}$ and 
$\mid\phi_{i}^{(m)}\mid^{2}: = 
\phi_{i}^{(m)}\cdot \overline{\phi_{i}^{(m)}}$. 

Then as in \cite{tu2}, we have the following lemma. 
\begin{lemma}(\cite[p.256, Theorem 1.1]{tu2}, see also \cite[Section 3.2]{tu3})
\[
K_{\infty}(B) := \mbox{the uppersemicontinuous envelope of} \,\,\, 
\overline{\lim}_{m\rightarrow\infty}\sqrt[m]{K_{m}(B)}
\]
exists and 
\[
h_{B,\infty} = \frac{1}{K_{\infty ,B}}. 
\]
 is an AZD of $K_{B}+X$.
\end{lemma}
We note that $K_{\infty}(B)$ may be totally different from 
$\overline{\lim}_{m\rightarrow\infty}\sqrt[m]{K_{m}(B)}$ on $X$, 
since we have taken the uppersemicontinuous envelope.
Let $\hat{X}$ denote the formal completion of $B$ along 
$X$. 
Then we see that 
\[
H^{0}(\hat{X},{\cal O}_{\hat{X}}(m(K_{B}+X)))
= \oplus_{\ell\geq 0}H^{0}(X,{\cal O}_{X}(mK_{X}+\ell A))
\]
holds.
We note that $h_{B,\infty}\mid_{X}$ and 
$h_{B,\infty}\mid_{X_{0}}$ are well defined, 
since $\kappa (X_{t}) \geq 0$ for every $t\in \Delta$. 
In fact by the assumption, there exists a positive integer
$m_{0}$ such that 
there exists a  nonzero section 
$\sigma_{0}\in H^{0}(X,{\cal O}_{X}(m_{0}K_{X}))$
which does not vanish identically on every fiber. 
Then $p^{*}\sigma_{0}$ is identified with 
a global holomorphic section of
${\cal O}_{B}(m_{0}(K_{B}+X))$ by the adjunction formula 
and the line bundle structure of $p : L\longrightarrow X$. 
Shrinking $\Delta$, if necessary, we may assume that 
$p^{*}\sigma_{0}$ is a bounded holomorphic section 
of $(m_{0}(K_{B}+X),h_{B})$ on $B$.
Let $B_{m}(1)$ denote the unit ball of ${\cal H}_{m}$. 
We note that 
\[
K_{m}(B)(z) = 
\sup_{\sigma\in B_{m}(1)} \mid \sigma\mid^{2}(z) \hspace{10mm} (z\in B)
\]
holds by definition. 
Then we see that there exists a positive constant $c_{0}$ 
such that 
\[
K_{\infty}(B) \geq c_{0}\cdot \mid p^{*}\sigma_{0}\mid^{2/m_{0}}
\]
holds on $B$. 
Hence $h_{B,\infty}\mid_{X}$ and $h_{B,\infty}\mid_{X_{0}}$
 are well defined. 

This implies that $h\mid_{X_{0}}$ is well defined by 
the definition of a canonical AZD. 
\begin{lemma}
\[
\bar{\cal I}(h^{m}) = \bar{\cal I}(h_{B,\infty}^{m}\mid_{X})
\]
and 
\[
\bar{\cal I}(h_{0}^{m}) = \bar{\cal I}(h^{m})\mid_{X_{0}} = \bar{\cal I}(h_{B,\infty}^{m}\mid_{X_{0}})
\]
hold for every $m\geq 0$.  
\end{lemma}
{\bf Proof of Lemma 4.2.}
We see that 
\[
\bar{\cal I}(h_{B,\infty}^{m}\mid_{X}) \subseteq 
\bar{\cal I}(h^{m})\mid_{X}
\]
holds for every $m\geq 0$ by the construction of $h$. 

On the other hand we see that 
\[
(\star )\hspace{10mm}\bar{\cal I}(h_{B,\infty}^{m}\mid_{X}) = \bar{\cal I}_{X}(h_{B,\infty}^{m})
\]
holds by Theorem 2.3, where 
\[
{\cal I}_{X}(h_{B,\infty}^{m}):= \cap_{\ell\geq 1}\sqrt[\ell]{{\cal I}(h_{B,\infty}^{m\ell})\mid_{D}}.
\]

Let $h_{m}$ be a canonical AZD of $mK_{X} +A$. 
Since $A$ is ample, we see that 
\[
\bar{\cal I}(h^{m\ell}) \subseteq \bar{\cal I}(h_{m\ell})
\]
holds for every $m,\ell\geq 0$. 
We note that for every $m,\ell \geq 0$, 
\[
{\cal O}_{X}(m\ell K_{X}+2A)\otimes {\cal I}(h_{m\ell})
\]
is generated by global sections by the definition of $A$.  
Hence shrinking $\Delta$ if necessary, we may and do assume that 
for every $m, \ell \geq 0$, 
\[
{\cal O}_{X}(m\ell K_{X}+2A)\otimes {\cal I}(h_{m\ell})
\]
is generated by 
\[
H^{0}_{(2)}(X,{\cal O}_{X}(m\ell K_{X}+2A)\otimes {\cal I}(h_{m\ell})).
\]
Hence noting the identity 
\[
H^{0}(\hat{X},{\cal O}_{\hat{X}}(m(K_{B}+X)))
= \oplus_{k\geq 0}H^{0}(X,{\cal O}_{X}(mK_{X}+kA))
\hspace{10mm} (m\geq 0),
\]
by the definition of $h_{B,\infty}$, we see that 
\[
\bar{\cal I}(h^{m}_{B,\infty})\mid_{D}  
\supseteq \overline{\cap_{\ell\geq 1}{\sqrt[\ell]{{\cal I}(h_{m\ell})}}}
\]
holds.
By the identity 
$(\star )$,  this implies that the opposite inclusion
\[
\bar{\cal I}(h^{m}) \subseteq \bar{\cal I}(h_{B,\infty}^{m}\mid_{X})
\]
holds for every $m\geq 0$.
Hence we conclude that the equality :
\[
\bar{\cal I}(h^{m}) = \bar{\cal I}(h_{B,\infty}^{m}\mid_{X})
\]
holds for every $m\geq 0$.

On the other hand, by Lemma 3.1 and $(\star )$, we see that 
for every $m$
\[
\bar{\cal I}(h_{B,\infty}^{m}\mid_{X_{0}})
\supseteq \overline{\cap_{\ell\geq 1}\sqrt[\ell]{{\cal I}(h_{\ell,0}^{m\ell})}}
\]
holds, where 
$h_{m,0}$ denotes a canonical AZD of $mK_{X_{0}}+A$
on $X_{0}$.
We note that obviously 
\[
\bar{\cal I}(h_{0}^{m\ell}) \subseteq \bar{\cal I}(h_{m,0}^{\ell})
\]
holds for every $m\geq 1$ and $\ell \geq 0$.
Hence we see that 
\[
\bar{\cal I}(h_{0}^{m}) \subseteq \bar{\cal I}(h_{B,\infty}^{m}\mid_{X_{0}})
\]
holds. 
On the other hand, since
\[
 \bar{\cal I}(h_{B,\infty}^{\ell})
= \bar{\cal I}(h^{\ell})
\]
holds for every $\ell\geq 0$ as we have seen above,
 by Theorem 2.3 we obtain that 
\[
 \bar{\cal I}(h_{B,\infty}^{m}\mid_{X_{0}})
= \bar{\cal I}(h^{m}\mid_{X_{0}})
\]
holds for every $m\geq 0$. 

Combining the above formulas, we see that 
\[
\bar{\cal I}(h_{0}^{m}) \subseteq  \bar{\cal I}(h_{B,\infty}^{m}\mid_{X_{0}})
= \bar{\cal I}(h^{m}\mid_{X_{0}})
\]
hold.
We note that by definition
\[
\bar{\cal I}(h^{m}\mid_{X_{0}})\subseteq \bar{\cal I}(h_{0}^{m}) 
\]
holds for every $m\geq 0$. 
Hence
\[
\bar{\cal I}(h_{0}^{m}) = \bar{\cal I}(h^{m}\mid_{X_{0}}) = \bar{\cal I}(h_{B,\infty}^{m}\mid_{X_{0}})
\]
hold for every $m\geq 0$. 
This completes the proof of Lemma 4.2. {\bf Q.E.D.} 

\subsection{Completion of the proof of Theorem 1.1}
By the $L^{2}$-extension theorem (\cite[p.200, Theorem]{o-t}), we see that 
\[
H^{0}(X,{\cal O}_{X}(mK_{X})\otimes {\cal I}(h_{B,\infty}^{m-1}\mid_{X}))
\rightarrow 
H^{0}(X_{0},{\cal O}_{X_{0}}(mK_{X_{0}})\otimes 
{\cal I}(h_{B,\infty}^{m}\mid_{X_{0}}))
\]
is surjective for every $m\geq 1$.

On the other hand, 
we note that 
\[
H^{0}(X_{0},{\cal O}_{X_{0}}(mK_{X_{0}}))
\simeq 
H^{0}(X_{0},{\cal O}_{X_{0}}(mK_{X_{0}})\otimes 
{\cal I}_{\infty}(h_{0}^{m}))
\]
holds for every $m\geq 0$ by the construction of 
the canonical AZD $h_{0}$ (cf. Remark 2.1). 
Since 
\[
{\cal I}_{\infty}(h_{0}^{m}) \subseteq 
\bar{\cal I}(h_{0}^{m}) \subseteq 
{\cal I}(h_{0}^{m})
\]
hold for every $m\geq 0$ and the closure of multiplier ideal 
sheaves are determined by the Lelong number of 
the pullback of the curvature current on some modifications as in Section 2.1, 
by Lemma 4.2
\[
H^{0}(X_{0},{\cal O}_{X_{0}}(mK_{X_{0}}))
\simeq 
H^{0}(X_{0},{\cal O}_{X_{0}}(mK_{X_{0}})\otimes 
{\cal I}(h_{B,\infty}^{m}\mid_{X_{0}}))
\]
holds for every $m\geq 0$. 
Hence we have that for every $m\geq 1$
\[
P_{m}(X_{0}) \leq P_{m}(X_{t})
\]
holds for every $t\in \Delta$. 
This means that $P_{m}(X_{t})$ is lowersemicontinuous
as a function on $t\in \Delta$. 
By the uppersemicontinuity theorem of cohomology groups, we see that 
$P_{m}(X_{t})$ is constant on $\Delta$. 
This completes the proof of Theorem 1.1.

Author's address\\
Hajime Tsuji\\
Department of Mathematics\\
Tokyo Institute of Technology\\
2-12-1 Ohokayama, Megro 152-8551\\
Japan \\
e-mail address: tsuji@math.titech.ac.jp
\end{document}